\documentclass{amsart}
\usepackage{amssymb}
\usepackage{mathrsfs}
\usepackage{enumerate}
\usepackage[cp1250]{inputenc}
\usepackage[OT4]{fontenc}
\newtheoremstyle{thm}{}{}{\it}{}{\sc}{.}{ }{}
\theoremstyle{thm}
\newtheorem{thm}{Theorem}
\newtheoremstyle{fact}{}{}{\it}{}{\sc}{.}{ }{}
\theoremstyle{fact}
\newtheorem{fact}{Fact}
\newtheoremstyle{rem}{}{}{}{}{\sc}{.}{ }{}
\theoremstyle{rem}
\newtheorem*{rem}{Remark}
\begin{document}
\title[Boundary behavior of the Kobayashi distance]{Boundary behavior of the Kobayashi distance in pseudoconvex Reinhardt domains}
\author{Tomasz Warszawski}
\address{Instytut Matematyki, Wydział Matematyki i Informatyki, Uniwersytet Jagielloński, ul. Prof. St. Łojasiewicza 6, 30-348 Kraków, Poland}
\email{tomasz.warszawski@im.uj.edu.pl}
\begin{abstract}
We prove that the Kobayashi distance near boundary of a pseudoconvex Reinhardt domain $D$ increases asymptotically at most like $-\log d_D+C$. Moreover, for boundary points from $\text{int}\overline D$ the growth does not exceed $\frac{1}{2}\log(-\log d_D)+C$. The lower estimate by $-\frac{1}{2}\log d_D+C$ is obtained under additional assumptions of $\mathcal C^1$-smoothness of a domain and a non-tangential convergence.
\end{abstract}
\maketitle
\section{Introduction and results}
The problem of a boundary behavior of the Kobayashi (pseudo)distance in pseudoconvex Reinhardt domains is connected with studying their Kobayashi completeness. The qualitative condition for the $k$-completeness of a bounded domain $D$ is $$k_{D}(z_0,z)\to\infty\text{ as } z\to\partial D.$$ The main fact is that if a pseudoconvex Reinhardt domain $D$ is hyperbolic then it is $k$-complete. At first Pflug \cite{pfl} proved it for bounded complete domains. A second step was done by Fu for bounded domains in \cite{fu}. The general case was finally solved by Zwonek in \cite{zw}. 

Hence it is natural to ask about a quantitative behavior of the function $k_{D}(z_0,\cdotp)$. Forstneri\v c and Rosay estimated it from below on bounded strongly pseudoconvex domains. Namely, it was proved in \cite{fr} that $$k_{D}(z_{1},z_{2})\geq-\frac{1}{2}\log d_{D}(z_{1})-\frac{1}{2}\log d_{D}(z_{2})+C$$ for $z_{j}$ near two distinct points $\zeta_{j}\in\partial D$, $j=1,2$. In the same paper the authors showed the opposite estimate for $\mathcal{C}^{1+\varepsilon}$-smooth domains with $z_{1},z_{2}$ near $\zeta_{0}\in\partial D$. This estimate in the bounded case follows from the inequality for the Lempert function of bounded $\mathcal{C}^{1+\varepsilon}$-smooth domains obtained by Nikolov, Pflug and Thomas $$\widetilde{k}_D(z_1,z_2)\leq-\frac{1}{2}\log d_D(z_1)-\frac{1}{2}\log d_D(z_2)+C,\,z_1,z_2\in D$$ in \cite{npt}. It was also proved that the above estimate fails in the $\mathcal{C}^1$-smooth case. The other general version of an upper estimate, for $\mathcal{C}^{2}$-smooth domains, can be found in \cite{jp}. The case of bounded convex domains was investigated by Mercer in \cite{mer}. For such domains we have $$-\frac{1}{2}\log d_{D}(z)+C'\leq k_{D}(z_{0},z)\leq-\alpha\log d_{D}(z)+C$$ with $\alpha>\frac{1}{2}$ and $z$ close to $\zeta_{0}\in\partial D$ (the constant $\alpha$ can not be replaced with $\frac{1}{2}$). An example $$D_\beta:=\{(z,w)\in\mathbb{C}^{2}:|z|^{\beta}+|w|^{\beta}<1\},\,0<\beta<1$$ shows that the lower estimate by $-\alpha\log d_{D}(z)+C$, where $\alpha>0$ --- a constant independent on a domain, is not true for complete pseudoconvex Reinhardt domains. Easy calculations lead to $$k_{D_\beta}((0,0),(z,0))\leq-\frac{\beta}{2}\log d_{D_\beta}(z,0)+C$$ if $0<z<1$ and $(z,0)$ tends to $(1,0)$.

In the paper we prove the following theorems.
\begin{thm}\label{1}
Let $D\subset\mathbb{C}^{n}$ be a pseudoconvex Reinhardt domain. Fix $z_{0}\in D$ and $\zeta_0\in\partial D$. Then for some constant $C$ the inequality $$k_{D}(z_{0},z)\leq-\log d_{D}(z)+C$$ holds if $z\in D$ tends to $\zeta_0$. 
Additionally, for $\zeta_0\in\mathbb{C}^{n}_{*}$ the estimate can be improved to $$k_{D}(z_{0},z)\leq-\frac{1}{2}\log d_{{D}}(z)+C'$$ where $C'$ is a constant.
\end{thm}
\begin{thm}\label{9}
Let ${D}\subset\mathbb{C}^{n}$ be a pseudoconvex Reinhardt domain. Fix $z_{0}\in D$ and $\zeta_{0}\in\partial D\cap\textnormal{int}\overline{D}$. Then for some constant $C$ the inequality $$k_{{D}}(z_{0},z)\leq\frac{1}{2}\log(-\log d_D(z))+C$$ holds if $z\in D$ tends to $\zeta_{0}$.
\end{thm}
\begin{thm}\label{3}
Let ${D}\subset\mathbb{C}^{n}$ be a $\mathcal{C}^{1}$-smooth pseudoconvex Reinhardt domain. Fix $z_{0}\in D$ and $\zeta_{0}\in\partial D$. Then for some constant $C$ the inequality $$k_{{D}}(z_{0},z)\geq-\frac{1}{2}\log d_{{D}}(z)+C$$ holds if $z\in D$ tends non-tangentially to $\zeta_{0}$.
\end{thm}
\section{Notations and definitions}
By $D$ we denote a domain in $\mathbb{C}^{n}$. The \textit{Kobayashi \emph{(}pseudo\emph{)}distance}\/ is defined as\medskip\\
$k_{D}(w,z):=\sup\{d_{D}(w,z):(d_{D})$ is a family of holomorphically invariant
\begin{flushright}
pseudodistances less than or equal to $\widetilde{k}_{D}\}$,\smallskip
\end{flushright}
where $$\widetilde{k}_{D}(w,z):=\inf\{p(\lambda,\mu):\lambda,\mu\in\mathbb{D}\ \textnormal{and }\exists f\in \mathcal{O}(\mathbb{D},D):\,f(\lambda)=w,\,f(\mu)=z\}$$
is the \textit{Lempert function}\/ of $D$, $\mathbb{D}\subset\mathbb{C}$ --- the unit disc and $p$ --- the Poincar\'e distance on $\mathbb{D}$. For general properties of functions $k_{D}$ one can see \cite{jp}. 

Denote $z_{j}$ as the $j$-th coordinate of point $z\in\mathbb{C}^{n}$. A domain $D$ is called a \textit{Reinhardt domain}\/ if $(\lambda_{1}z_{1},\ldots,\lambda_{n}z_{n})\in D$ for all numbers $\lambda_{1},\ldots,\lambda_{n}\in\partial\mathbb{D}$ and points $z\in D$. A Reinhardt domain $D$ is \textit{complete in $j$-th direction}\/ if $$(\{1\}^{j-1}\times\overline{\mathbb{D}}\times\{1\}^{n-j})\,\cdotp D\subset D,$$ where $A\,\cdotp B:=\{(a_{1}b_{1},\ldots,a_{n}b_{n}):a\in A,\,b\in B\}$. Define subspaces $V_{j}^{n}:=\{z\in\mathbb{C}^{n}:z_{j}=0\}$ for $j=1,\ldots,n$. If a Reinhardt domain $D$ is complete in the $j$-th direction for all $j$ such that $D\cap V_{j}^{n}\neq\emptyset$ then $D$ is called \textit{relatively complete}.

Let us denote $A_{*}:=A\setminus\{0\}$ for a set $A\subset\mathbb{C}$ and $\mathbb{C}^{n}_{*}:=(\mathbb{C}_{*})^{n}$. By $d_{D}(z)$ denote a distance of a point $z\in D$ to $\partial D$ (here, exceptionally, $D$ can be a domain in $\mathbb{R}^{n}$) and by $\zeta_{D}(z)$ --- one of points admitting a distance of a point $z\in D$ to $\partial D$. 

We will use the following main branch of the power $z^{\alpha}:=e^{\alpha\log z}=e^{\alpha(\log|z|+i\textnormal{Arg}\,z)}$, where the main argument $\textnormal{Arg}\,z\in(-\pi,\pi]$. Define $z^{\alpha}:=z_{1}^{\alpha_{1}}\cdotp\ldots\cdotp z_{n}^{\alpha_{n}}$, $|z|^{\alpha}:=|z_{1}|^{\alpha_{1}}\cdotp\ldots\cdotp |z_{n}|^{\alpha_{n}}$ for $z\in\mathbb{C}^{n}_{*}$ and $\alpha\in\mathbb{R}^{n}$. Moreover, let $|z|:=(|z_{1}|,\ldots,|z_{n}|)$ for $z\in\mathbb{C}^{n}$, $\log|z|:=(\log|z_{1}|,\ldots,\log|z_{n}|)$ for $z\in\mathbb{C}^{n}_{*}$ and $\log D:=\{\log|z|:z\in D\cap\mathbb{C}^{n}_{*}\}$ --- a \textit{logarithmic image}\/ of $D$. We use $C$ to denote constants not necessarily the same in different places. We also need notations $f\lesssim g$ if there exists $C>0$ such that $f\leq Cg$; $f\approx g$ if $f\lesssim g$ and $g\lesssim f$. 

We call $D$ a $\mathcal{C}^{k}$\textit{-smooth}\/ domain if for any point $\zeta_{0}\in\partial D$ there exist its open neighbourhood $U\subset\mathbb{C}^{n}$ and a $\mathcal{C}^{k}$-smooth function $\rho:U\longrightarrow\mathbb{R}$ such that
\begin{enumerate}
\item $U\cap D=\{z\in U:\rho(z)<0\}$;
\item $U\setminus\overline{D}=\{z\in U:\rho(z)>0\}\neq\emptyset$;
\item $\nabla\rho:=\left(\frac{\partial\rho}{\partial\overline{z}_{1}},\ldots,\frac{\partial\rho}{\partial \overline{z}_{n}}\right)\neq 0$ on $U$.
\end{enumerate}
The function $\rho$ is called a \textit{local defining function}\/ for $D$ at the point $\zeta_{0}$. 

For a $\mathcal{C}^{1}$-smooth domain $D$ we define a \textit{normal vector}\/ to $\partial D$ at a point $\zeta_{0}\in\partial D$ as $$\nu_{D}(\zeta_{0}):=\frac{\nabla\rho(\zeta_{0})}{\|\nabla\rho(\zeta_{0})\|},$$
where $\rho$ is a local defining function for $D$ at $\zeta_{0}$. Clearly $$z=\zeta_{D}(z)-d_{D}(z)\nu_{D}(\zeta_{D}(z))$$ for $z\in D$ and
$$\lim_{D\ni z\to\zeta_{0}}\nu_{D}(\zeta_{D}(z))=\nu_{D}(\zeta_{0})$$ for every choice of $\zeta_{D}(z)$. For the transparent notation we shorten the symbol $\nu_{D}(\zeta_{D}(z))$ to $\nu_{D}(z)$.

To define a non-tangential convergence we need a concept of a $\textit{cone}$ with a vertex $x_{0}\in\mathbb{R}^{n}$, a semi-axis $\nu\in(\mathbb{R}^{n})_{*}$ and an angle $\alpha\in (0,\frac{\pi}{2})$. It is a set of $x\in\mathbb{R}^{n}\setminus\{x_{0}\}$ such that an angle between vectors $\nu$ and $x-x_{0}$ does not exceed $\alpha$. Let $D$ be a $\mathcal{C}^{1}$-smooth domain and $\zeta_{0}\in\partial D$. We say that $z\in D$ tends \textit{non-tangentially}\/ to $\zeta_{0}$ if there exist a cone $\mathcal{A}\subset\mathbb{C}^{n}\cong\mathbb{R}^{2n}$ with a vertex $\zeta_{0}$, a semi-axis $-\nu_{D}(\zeta_{0})$ and an angle $\alpha\in (0,\frac{\pi}{2})$ and an open neighbourhood $U\subset\mathbb{C}^{n}$ of $\zeta_{0}$ such that $U\cap\mathcal{A}\subset D$ and $z$ tends to $\zeta_{0}$ in $U\cap\mathcal{A}$.

We say that a Reinhardt domain $D$ satisfies the \textit{Fu condition}\/ if for any $j\in\{1,\ldots,n\}$ the following implication holds $$\partial D\cap V_{j}^{n}\neq\emptyset\Longrightarrow D\cap V_{j}^{n}\neq\emptyset.$$

The following well-known properties of pseudoconvex Reinhardt domains will be used in the paper (see e.g. \cite{jp1}). 
\begin{fact}\label{4}
A Reinhardt domain $D$ is pseudoconvex if and only if $\,\log D$ is convex and $D$ relatively complete.
\end{fact}
\begin{fact}\label{5}
A $\mathcal{C}^{1}$-smooth Reinhardt domain satisfies the Fu condition.
\end{fact}
\section{Proofs}
\begin{proof}[\sc Proof of Theorem \ref{1}] 
We proceed as follows. The first step is to simplify the general case to `real' coordinates, further we consider some parallelepipeds contained in the given domain and use the decreasing property of the Kobayashi distance. Finally, we explicitly calculate and estimate it in other domains --- cartesian products of a belt and annuli in $\mathbb{C}$. To improve the estimate for a boundary point with all non-zero coordinates we use similar methods, but with intervals instead of parallelepipeds.

Using some biholomorphism of the form $$w\ni\mathbb{C}^{n}\longmapsto(a_{1}w_{1},\ldots,a_{n}w_{n})\in\mathbb{C}^{n},\,a\in\mathbb{C}_{*}^{n}$$ and the triangle inequality for $k_{D}$, we can assume that $z_{0}=(1,\ldots,1)$ and $|\zeta_{0j}|\neq 1$ for $j=1,\ldots,n$. Notice that the proof can be reduced to $z\in D\cap\mathbb{C}_{*}^{n}$ near $\zeta_0$ and next to the case $$z\in D\cap(0,\infty)^n \text{ near } \zeta_0\in\partial D\cap([0,\infty)\setminus\{1\})^n.$$ Indeed, the first reduction follows from the continuity of $k_{D}$ and the triangle inequality for $k_{D}$. Now, if $z\to\zeta_0$ then $|z|\to|\zeta_0|\in\partial D$ and $$k_D(z_0,z)=k_D(\widetilde{z_0},|z|),$$ where
$$\widetilde{z_0}:=\left(\frac{|z_1|}{z_1}z_{01},\ldots,\frac{|z_n|}{z_n}z_{0n}\right)\in T:=\{(\lambda_{1}z_{01},\ldots,\lambda_{n}z_{0n}):\lambda_{1},\ldots,\lambda_{n}\in\partial\mathbb{D}\}.$$ The continuity of $k_{D}$ gives $$\max_{T\times T}k_{D}=:C<\infty$$ and therefore $$k_D(\widetilde{z_0},|z|)\leq k_D(\widetilde{z_0},z_0)+k_D(z_0,|z|)\leq k_D(z_0,|z|)+C.$$
The property $d_D(|z|)=d_D(z)$ finishes this reduction. In what follows, we assume that points $z\in D\cap(0,\infty)^n$ are sufficiently close to $\zeta_0\in\partial D\cap([0,\infty)\setminus\{1\})^n$.

Observe that $$d_{\log D}(\log z)\geq\varepsilon d_{D}(z)$$ for some $\varepsilon>0$. Indeed, for $u\in\mathbb{R}^{n} $, $\|u\|<1$ and $0\leq t\leq\varepsilon d_{D}(z)$, where $$\varepsilon:=\frac{1}{3(\|\zeta_0\|+1)},$$ we have
$$\log z+tu\in\log D$$ if and only if $$\left(z_{1}e^{tu_{1}},\ldots,z_{n}e^{tu_{n}}\right)\in D$$ but this property follows from $$\left\|(z_{j}e^{tu_{j}})_{j=1}^{n}-z\right\|\leq\sqrt{\sum_{j=1}^{n}z_{j}^{2}(2t)^{2}}\leq 2t(\|\zeta_0\|+1)<d_{D}(z).$$
Moreover, for $\zeta_0=0$ a similar consideration leads to $$d_{\log D}(\log z)\geq\varepsilon'\frac{d_{D}(z)}{\|z\|}$$ for sufficiently small $\varepsilon'>0$. Indeed, there exists $\varepsilon'\in\left(0,\frac{1}{2}\right)$ such that the inequalities $$\left|e^{tu_{j}}-1\right|\leq 2t,\,j=1,\ldots,n$$ hold for $0\leq t\leq\varepsilon'$. Hence for $0\leq t\leq\varepsilon'\frac{d_{D}(z)}{\|z\|}$ we have $$\left\|(z_{j}e^{tu_{j}})_{j=1}^{n}-z\right\|\leq\sqrt{\sum_{j=1}^{n}z_{j}^{2}(2t)^{2}}\leq 2\varepsilon'\frac{d_{D}(z)}{\|z\|}\|z\|<d_{D}(z).$$

Denote $$\widetilde{d}_{D}(z):=\begin{cases}\varepsilon d_{D}(z),\,\zeta_0\neq 0\\\varepsilon''\frac{d_{D}(z)}{\|z\|},\,\zeta_0=0,\end{cases}$$ where $\varepsilon'':=\varepsilon' d_{\log D}(0)$. 
Let us define $$m_{z}:=\min\{0,\log z_{1}\},\,M_{z}:=\max\{0,\log z_{1}\}$$ and consider the set $$D_{z}:=\{w\in\mathbb{C}^{n}:m_{z}-\widetilde{d}_{D}(z)<\log|w_{1}|<M_{z}+\widetilde{d}_{D}(z),$$$$\frac{\log z_{j}}{\log z_{1}}\log|w_{1}|-\widetilde{d}_{D}(z)<\log|w_{j}|<\frac{\log z_{j}}{\log z_{1}}\log|w_{1}|+\widetilde{d}_{D}(z),\,j=2,\ldots,n\}.$$
Then $\log D_{z}$ is a domain in $\mathbb{R}^{n}$ containing points $0$ and $\log z$ but contained in a convex domain $\log D$.
Define also $$G_{z}:=\{v\in\mathbb{C}^{n}:m_{z}-\widetilde{d}_{D}(z)<\textnormal{Re\,}v_{1}<M_{z}+\widetilde{d}_{D}(z),
$$$$-\widetilde{d}_{D}(z)<\log|v_{j}|<\widetilde{d}_{D}(z),\,j=2,\ldots,n\}.$$
Hence the holomorphic map $$f_{z}(v):=\left(e^{v_{1}},v_{2}{e^{v_{1}\frac{\log z_{2}}{\log z_{1}}}},\ldots,
v_{n}{e^{v_{1}\frac{\log z_{n}}{\log z_{1}}}}\right),\,v\in G_{z}$$
has values in $D_{z}$. Moreover $$w=f_{z}\left(\log w_{1},\frac{w_{2}}{w_{1}^{\frac{\log z_{2}}{\log z_{1}}}},\ldots,\frac{w_{n}}{w_{1}^{\frac{\log z_{n}}{\log z_{1}}}}\right)\textnormal{ for \,} w\in D_{z}.$$ Therefore $$k_{D}(z_{0},z)\leq k_{D_{z}}(z_{0},z)=k_{D_{z}}\left(f_{z}(0,1,\ldots,1),f_{z}\left(\log z_{1},\frac{z_{2}}{z_{1}^{\frac{\log z_{2}}{\log z_{1}}}},\ldots,\frac{z_{n}}{z_{1}^{\frac{\log z_{n}}{\log z_{1}}}}\right)\right)=$$$$=k_{D_{z}}(f_{z}(0,1,\ldots,1),f_{z}(\log z_{1},1,\ldots,1))\leq k_{G_{z}}((0,1,\ldots,1),(\log z_{1},1,\ldots,1))=$$
$$=\max\{k_{S_{z}}(0,\log z_{1}),k_{A_{z}}(1,1),\ldots,
k_{A_{z}}(1,1)\}=k_{S_{z}}(0,\log z_{1}),$$ where $$S_{z}:=\left\{\lambda\in\mathbb{C}:m_{z}-\widetilde{d}_{D}(z)<\textnormal{Re\,}\lambda<M_{z}+\widetilde{d}_{D}(z)\right\}$$
and $$A_{z}:=\left\{\lambda\in\mathbb{C}:-\widetilde{d}_{D}(z)<\log|\lambda|<\widetilde{d}_{D}(z)\right\}.$$
Using suitable biholomorphisms, we calculate 
$$k_{S_{z}}(0,\log z_{1})=p\left(\frac{i-\exp\pi iP(z)}{i+\exp\pi iP(z)},\frac{i-\exp\pi iQ(z)}{i+\exp\pi iQ(z)}\right),$$ where $$P(z):=\frac{\widetilde{d}_{D}(z)-m_{z}}{2\widetilde{d}_{D}(z)+M_{z}-m_{z}},\,Q(z):=\frac{\log z_{1}+\widetilde{d}_{D}(z)-m_{z}}{2\widetilde{d}_{D}(z)+M_{z}-m_{z}}.$$

Analogously, after changing the index 1 to any of $2,\ldots,n$, we get $$k_{D}(z_{0},z)\leq\min_{j=1,\ldots,n}k_{S_{z}^{(j)}}(0,\log z_{j}),$$ where $$k_{S_{z}^{(j)}}(0,\log z_{j})=p\left(\frac{i-\exp\pi iP^{(j)}(z)}{i+\exp\pi iP^{(j)}(z)},\frac{i-\exp\pi iQ^{(j)}(z)}{i+\exp\pi iQ^{(j)}(z)}\right)$$ and $$S_{z}^{(j)}:=\left\{\lambda\in\mathbb{C}:m_{z}^{(j)}-\widetilde{d}_{D}(z)<\textnormal{Re\,}
\lambda<M_{z}^{(j)}+\widetilde{d}_{D}(z)\right\},$$
$$m_{z}^{(j)}:=\min\{0,\log z_{j}\},\,M_{z}^{(j)}:=\max\{0,\log z_{j}\},\,j=1,\ldots,n,$$ 
$$P^{(j)}(z):=\frac{\widetilde{d}_{D}(z)-m_{z}^{(j)}}{2\widetilde{d}_{D}(z)+M_{z}^{(j)}-m_{z}^{(j)}},\,Q^{(j)}(z):=\frac{\log z_{1}+\widetilde{d}_{D}(z)-m_{z}^{(j)}}{2\widetilde{d}_{D}(z)+M_{z}^{(j)}-m_{z}^{(j)}}.$$

Consider two cases: $\zeta_0\neq 0$ and $\zeta_0=0$. If $\zeta_0\neq 0$ then choose $j\in\{1,\ldots,n\}$ such that $\zeta_{0j}\neq 0$ (recall that $\zeta_{0j}=|\zeta_{0j}|\neq 1$). In the case of $\zeta_{0j}>1$ we obtain \begin{equation}\label{6}k_{S_{z}^{(j)}}(0,\log z_{j})=p\left(\frac{i-\exp\pi iT^{(j)}(z)}{i+\exp\pi iT^{(j)}(z)},\frac{i-\exp\pi iU^{(j)}(z)}{i+\exp\pi iU^{(j)}(z)}\right)\leq
\end{equation}
$$\leq p\left(0,\frac{i-\exp\pi iT^{(j)}(z)}{i+\exp\pi iT^{(j)}(z)}\right)+p\left(0,\frac{i-\exp\pi iU^{(j)}(z)}{i+\exp\pi iU^{(j)}(z)}\right),$$ where $$T^{(j)}(z):=\frac{\varepsilon d_{D}(z)}{2\varepsilon d_{D}(z)+\log z_{j}},\,U^{(j)}(z):=\frac{\log z_{j}+\varepsilon d_{D}(z)}{2\varepsilon d_{D}(z)+\log z_{j}}.$$  We have, by Taylor expansion
$$\frac{i-\exp\pi iT^{(j)}(z)}{i+\exp\pi iT^{(j)}(z)}=i-\pi iT^{(j)}(z)+O\left(d_D(z)^2\right).$$ Hence 
$$p\left(0,\frac{i-\exp\pi iT^{(j)}(z)}{i+\exp\pi iT^{(j)}(z)}\right)=p\left(0,i-\pi iT^{(j)}(z)+O\left(d_D(z)^2\right)\right)\leq$$$$\leq\frac{\log 2}{2}-\frac{1}{2}\log\left(1-\left|i-\pi iT^{(j)}(z)+O\left(d_D(z)^2\right)\right|\right)\leq$$$$\leq\frac{\log 2}{2}-\frac{1}{2}\log\left(1-\left|i-\pi iT^{(j)}(z)\right|-\left|O\left(d_D(z)^2\right)\right|\right)=$$$$=\frac{\log 2}{2}-\frac{1}{2}\log\left(\pi\frac{\varepsilon d_{D}(z)}{2\varepsilon d_{D}(z)+\log z_{j}}-O\left(d_D(z)^2\right)\right)\leq-\frac{1}{2}\log d_{D}(z)+C.$$ Similarly $$\frac{i-\exp\pi iU^{(j)}(z)}{i+\exp\pi iU^{(j)}(z)}=-i+\pi iT^{(j)}(z)+O\left(d_D(z)^2\right),$$ which gives the same estimation for the second summand.

Otherwise if $\zeta_{0j}<1$, we have
\begin{equation}\label{7}k_{S_{z}^{(j)}}(0,\log z_{j})=p\left(\frac{i-\exp\pi iV^{(j)}(z)}{i+\exp\pi iV^{(j)}(z)},\frac{i-\exp\pi iW^{(j)}(z)}{i+\exp\pi iW^{(j)}(z)}\right),
\end{equation} where $$V^{(j)}(z):=\frac{\varepsilon d_{D}(z)-\log z_{j}}{2\varepsilon d_{D}(z)-\log z_{j}},\,W^{(j)}(z):=\frac{\varepsilon d_{D}(z)}{2\varepsilon d_{D}(z)-\log z_{j}}.$$
We see that the expression in \eqref{7} is the expression in \eqref{6} after substitute $\log z_j\leadsto-\log z_j$ and the estimates stay true.

Assume $\zeta_0=0$. We have for $j=1,\ldots,n$
$$k_{S_{z}^{(j)}}(0,\log z_{j})=p\left(\frac{i-\exp\pi iX^{(j)}(z)}{i+\exp\pi iX^{(j)}(z)},\frac{i-\exp\pi iY^{(j)}(z)}{i+\exp\pi iY^{(j)}(z)}\right)\leq$$$$\leq p\left(0,\frac{i-\exp\pi iX^{(j)}(z)}{i+\exp\pi iX^{(j)}(z)}\right)+p\left(0,\frac{i-\exp\pi iY^{(j)}(z)}{i+\exp\pi iY^{(j)}(z)}\right),$$ where
$$X^{(j)}(z):=\frac{\varepsilon'' d_{D}(z)\|z\|^{-1}-\log z_{j}}{2\varepsilon'' d_{D}(z)\|z\|^{-1}-\log z_{j}},\,Y^{(j)}(z):=\frac{\varepsilon'' d_{D}(z)\|z\|^{-1}}{2\varepsilon'' d_{D}(z)\|z\|^{-1}-\log z_{j}}.$$
Putting $$\delta^{(j)}(z):=\frac{\varepsilon'' d_{D}(z)}{\|z\|\log z_{j}}$$ we have $$X^{(j)}(z)=\frac{\delta^{(j)}(z)-1}{2\delta^{(j)}(z)-1},\,Y^{(j)}(z)=\frac{\delta^{(j)}(z)}{2\delta^{(j)}(z)-1}$$ and $\delta^{(j)}(z)\to 0$ as $z\to 0$. The analogous calculations as in the first case give $$\frac{i-\exp\pi iX^{(j)}(z)}{i+\exp\pi iX^{(j)}(z)}=-i+\pi iY^{(j)}(z)+O\left(\delta^{(j)}(z)^2\right)$$ and $$\frac{i-\exp\pi iY^{(j)}(z)}{i+\exp\pi iY^{(j)}(z)}=i-\pi iY^{(j)}(z)+O\left(\delta^{(j)}(z)^2\right).$$ Therefore 
$$p\left(0,\frac{i-\exp\pi iX^{(j)}(z)}{i+\exp\pi iX^{(j)}(z)}\right)\leq\frac{\log 2}{2}-\frac{1}{2}\log\left(\pi\frac{\delta^{(j)}(z)}{2\delta^{(j)}(z)-1}-O\left(\delta^{(j)}(z)^2\right)\right)\leq$$
$$\leq-\frac{1}{2}\log(-\delta^{(j)}(z))+C$$ and similarly $$p\left(0,\frac{i-\exp\pi iY^{(j)}(z)}{i+\exp\pi iY^{(j)}(z)}\right)\leq-\frac{1}{2}\log(-\delta^{(j)}(z))+C.$$ Finally $$\min_{j=1,\ldots,n}k_{S_{z}^{(j)}}(0,\log z_{j})\leq\min_{j=1,\ldots,n}-\log(-\delta^{(j)}(z))+C=$$$$=-\log d_{D}(z)+\log\|z\|+\min_{j=1,\ldots,n}\log(-\log z_{j})+C=$$$$=-\log d_{D}(z)+\log\|z\|+\log\left(-\log\max_{j=1,\ldots,n}z_{j}\right)+C\leq$$$$\leq-\log d_{D}(z)+\log\|z\|+\log(-\log\|z\|)+C\leq-\log d_D(z)+C.$$ 

For improving the estimate in the case of $\zeta_0\in\partial D\cap\mathbb{C}^{n}_{*}$, we may assume that $z_0\in\mathbb{C}_{*}^{n}$ and $|z_{0j}|,|\zeta_{0j}|\neq 1$ for $j=1,\ldots,n$. Since $\log D$ is a convex domain, the interval $$I_{z}:=\{t\log|z|+(1-t)\log|z_0|:t\in(-\varepsilon(z),1+\delta(z))\}$$ is contained in $\log D$ for some positive numbers $\delta(z)$, $\varepsilon(z)$. The number $\varepsilon(z)$ can be chosen as a sufficiently small positive constant $\varepsilon$ independent of $z$. Indeed, 
$$t\log|z|+(1-t)\log|z_0|=\log|z_0|+t(\log|z|-\log|z_0|)$$ and $\|\log|z|-\log|z_0|\|$ is bounded, say by $M$. Hence $$\varepsilon:=\frac{d_{\log D}(\log|z_0|)}{2M}$$ is good. Analogously, $$\frac{d_{\log D}(\log|z|)}{2M}$$ is a candidate for $\delta(z)$. We have $$\frac{d_{\log D}(\log|z|)}{2M}\geq\delta d_{D}(z)$$ for some $\delta>0$ (in fact, ``$\geq$'' can be replaced with ``$\approx$''). Thus we can choose $\delta(z):=\delta d_{D}(z)$. 

From the inclusion $I_{z}\subset\log D$ it follows that $$\exp I_{z}\subset D$$ i.e.
$$\left(\left|\frac{z_{1}}{z_{01}}\right|^{t}|z_{01}|,\ldots,\left|\frac{z_{n}}{z_{0n}}\right|^{t}|z_{0n}|\right)\in D$$ for $t\in(-\varepsilon,1+\delta d_{D}(z))$. Hence the holomorphic map
$$f_z(\lambda):=\left(e^{i\arg z_{1}}\left|\frac{z_{1}}{z_{01}}\right|^{\lambda}|z_{01}|,\ldots,e^{i\arg z_{n}}\left|\frac{z_{n}}{z_{0n}}\right|^{\lambda}|z_{0n}|\right)$$ leading from the strip $$S_{z}:=\{\lambda\in\mathbb{C}:-\varepsilon<\textnormal{Re}\,\lambda<1+\delta d_{D}(z)\}$$ has values in $D$. Moreover $f_z(1)=z$ and $f_z(0)$ lies on the torus $$T:=\{(\lambda_{1}z_{01},\ldots,\lambda_{n}z_{0n}):\lambda_{1},\ldots,\lambda_{n}\in\partial\mathbb{D}\}.$$
Therefore $$k_{D}(z_{0},z)\leq k_{D}(z_{0},f_z(0))+k_{D}(f_z(0),z)\leq$$$$\leq k_{D}(f_z(0),f_z(1))+\max_{T\times T}k_{D}\leq k_{S_{z}}(0,1)+\max_{T\times T}k_{D}.$$
Calculating $k_{S_{z}}(0,1)$ we get $$k_{S_{z}}(0,1)=p\left(\frac{i-\exp\pi iP^{(j)}(z)}{i+\exp\pi iP^{(j)}(z)},\frac{i-\exp\pi iQ^{(j)}(z)}{i+\exp\pi iQ^{(j)}(z)}\right),$$
where $$P^{(j)}(z):=\frac{\varepsilon}{1+\varepsilon+\delta d_{D}(z)},\,Q^{(j)}(z):=\frac{1+\varepsilon}{1+\varepsilon+\delta d_{D}(z)}.$$

Certainly, first of the above argument of the function $p$ tends to some point from the unit disc. For the second we have $$\frac{i-\exp\pi iQ^{(j)}(z)}{i+\exp\pi iQ^{(j)}(z)}=-i+\pi i\frac{\delta d_{D}(z)}{1+\varepsilon+\delta d_{D}(z)}+O\left(d_D(z)^2\right).$$ Consequently $$p\left(0,\frac{i-\exp\pi iQ^{(j)}(z)}{i+\exp\pi iQ^{(j)}(z)}\right)\leq\frac{\log 2}{2}-\frac{1}{2}\log\left(\pi\frac{\delta d_{D}(z)}{1+\varepsilon+\delta d_{D}(z)}-O\left(d_{D}(z)^2\right)\right)\leq$$$$\leq-\frac{1}{2}\log d_{D}(z)+C.$$ 
The triangle inequality for $p$ finishes the proof. 
\end{proof}
\begin{proof}[\sc Proof of Theorem \ref{9}] 
The proof is based on decreasing and product properties of the Kobayashi distance and need to consider some cases which form, in fact, an induction.

Note that if $E\subset\mathbb{R}^n$ is a convex domain then $E=\text{int}\overline{E}$. The condition $\zeta_{0}\in\partial D\cap\text{int}\overline{D}$ implies $\zeta_{0}\notin\mathbb{C}^{n}_{*}$. To see this, assume that $\zeta_{0}\in\mathbb{C}^{n}_{*}$. An easy topological argument shows that $$\log|\zeta_0|\in(\partial\log D)\cap\text{int\,}\overline{\log D}=(\partial\log D)\cap\log D=\emptyset.$$

Assume, without loss of generality, that $$\zeta_{0}=(\zeta_{01},\ldots,\zeta_{0k},0,\ldots,0),$$ where $0\leq k\leq n-1$ and $\zeta_{0j}\neq 0$, $j\leq k$. Let $r>0$ be such that an open polydisc $P(\zeta_0,r)$ is contained in $\overline{D}$. Then $\log P(\zeta_0,r)\subset\log\overline{D}$. Taking interiors of both sides we get $$\log P(\zeta_0,r)\subset\text{int}\log\overline{D}=\text{int\,}\overline{\log D}=\log D.$$ Therefore
\begin{equation}\label{10}P(\zeta_0,r)\cap\mathbb{C}^{n}_{*}\subset D.\end{equation} Clearly (for fixed small $r$)$$P(\zeta_0,r)\cap\mathbb{C}^{n}_{*}=\mathbb{D}(\zeta_{01},r)\times\ldots\times\mathbb{D}(\zeta_{0k},r)\times (r\mathbb{D}_*)^{n-k},$$ where $\mathbb{D}(\zeta_{0j},r)$ is a disc in $\mathbb{C}$ centered at $\zeta_{0j}$ with radius $r$. Hence, choosing any $z_0\in P(\zeta_0,r)\cap\mathbb{C}^{n}_{*}$, we have $$k_D(z_0,z)\leq\max\left\{\max_{j=1,\ldots,k}k_{\mathbb{D}(\zeta_{0j},r)}(z_{0j},z_j),\max_{j=k+1,\ldots,n}
k_{r\mathbb{D}_*}(z_{0j},z_{j})\right\}$$ for $z\in D\cap\mathbb{C}^{n}_{*}$ near $\zeta_0$. For $j=1,\ldots,k$ the numbers $z_j$ tend to $\zeta_{0j}$, so the first of the above maxima is bounded by a constant. The well-known estimate for the punctured disc gives us $$k_{r\mathbb{D}_*}(z_{0j},z_{j})\leq\frac{1}{2}\log(-\log d_{r\mathbb{D}_*}(z_j))+C=
\frac{1}{2}\log(-\log|z_j|)+C$$ for $j=k+1,\ldots,n$. Therefore 
\begin{equation}\label{11}k_D(z_0,z)\leq\frac{1}{2}\log\left(-\log\min_{j=k+1,\ldots,n}|z_j|\right)+C.\end{equation}

The above estimate is not sufficiently good yet. Denote $z':=(z_1,\ldots,z_k)$. Note that \begin{equation}\label{13}(z',0,\ldots,0)\in\partial D.\end{equation} Indeed, $(z',0,\ldots,0)\in\overline{D}$. If $(z',0,\ldots,0)\in D$ then $D$ is complete in the directions $k+1,\ldots,n$ (Fact \ref{4}). Moreover, $(\zeta_{01},\ldots,\zeta_{0k},r/2,\ldots,r/2)\in D$, which implies $(\zeta_{01},\ldots,\zeta_{0k},0,\ldots,0)\in D$ --- a contradiction.

We claim that for all $k+1\leq p<q\leq n$ \begin{equation}\label{12}(z',0,\ldots,0,\underline{z_p},0,\ldots,0)\in\partial D \text{ or }(z',0,\ldots,0,\underline{z_q},0,\ldots,0)\in\partial D,\end{equation}
where the symbol $\underline{z_j}$ means that $z_j$ is on the $j$-th place. If \eqref{12} it is not true then both points belong to $D$ (recall that $P(\zeta_0,r)\subset\overline{D}$). Hence $D$ is complete in the directions $k+1,\ldots,n$ and $(z',0,\ldots,0)\in D$, which contradicts \eqref{13}. 

Therefore all points $$(z',0,\ldots,0,\underline{z_p},0,\ldots,0),\,p=k+1,\ldots,n,$$ except possibly one, belong to $\partial D$. Consider the following cases.

Case 1.1. One of above points, say $(z',0,\ldots,0,z_n)$, does not belong to $\partial D$. Then it belongs to $D$. Hence $D$ is complete in the directions $k+1,\ldots,n-1$. Now the inclusion \eqref{10} can be improved to $$P(\zeta_0,r)\cap(\mathbb{C}^{n-1}\times\mathbb{C}_*)\subset D$$ and $$P(\zeta_0,r)\cap(\mathbb{C}^{n-1}\times\mathbb{C}_*)=\mathbb{D}(\zeta_{01},r)\times\ldots\times\mathbb{D}(\zeta_{0k},r)\times (r\mathbb{D})^{n-k-1}\times r\mathbb{D}_*.$$ The estimate for $k_D(z_0,z)$ is improved to $$\max\left\{\max_{j=1,\ldots,k}k_{\mathbb{D}(\zeta_{0j},r)}(z_{0j},z_j),\max_{j=k+1,\ldots,n-1}
k_{r\mathbb{D}}(z_{0j},z_{j}),\,k_{r\mathbb{D}_*}(z_{0n},z_{n})\right\}=$$$$= k_{r\mathbb{D}_*}(z_{0n},z_{n})\leq\frac{1}{2}\log(-\log|z_n|)+C.$$ It remains to notice that $$(z',z_{k+1},\ldots,z_{n-1},0)\in\partial D$$ since in the opposite case the domain $D$ would be complete in $n$-th direction and the property $(z',0,\ldots,0,z_n)\in D$ would imply $(z',0,\ldots,0)\in D$ --- a contradiction with \eqref{13}. Thus $$d_D(z)\leq\|z-(z',z_{k+1},\ldots,z_{n-1},0)\|=|z_n|,$$ which let us estimate $$\frac{1}{2}\log(-\log|z_n|)+C\leq\frac{1}{2}\log(-\log d_D(z))+C.$$

Case 1.2. All the points $$(z',0,\ldots,0,\underline{z_p},0,\ldots,0),\,p=k+1,\ldots,n$$ belong to $\partial D$. We claim that for all $k+1\leq p<q\leq n$ and $k+1\leq p'<q'\leq n$ with $\{p,q\}\neq\{p',q'\}$ $$(z',0,\ldots,0,\underline{z_p},0,\ldots,0,\underline{z_q},0,\ldots,0)\in\partial D \text{ or }(z',0,\ldots,0,\underline{z_{p'}},0,\ldots,0,\underline{z_{q'}},0,\ldots,0)\in\partial D.$$ Analogously as before we use an argument of completeness in the suitable directions to get $$(z',0,\ldots,0,\underline{z_j},0,\ldots,0)\in D$$ for some $j\in\{p,q,p',q'\}$ --- a contradiction with the assumption of the case 1.2. Therefore all points $$(z',0,\ldots,0,\underline{z_p},0,\ldots,0,\underline{z_q},0,\ldots,0),\,k+1\leq p<q\leq n,$$ except possibly one, belong to $\partial D$. Again we consider two cases.

Case 2.1. One of above points, say $(z',0,\ldots,0,z_{n-1},z_n)$, does not belong to $\partial D$. Then it belongs to $D$. We see, analogously as in the case 1.1, that $$P(\zeta_0,r)\cap(\mathbb{C}^{n-2}\times\mathbb{C}^2_*)\subset D,$$ $$k_D(z_0,z)\leq\frac{1}{2}\log\left(-\log\min_{j=n-1,n}|z_j|\right)+C,$$ $$(z',z_{k+1},\ldots,z_{n-2},z_{n-1},0),\,(z',z_{k+1},\ldots,z_{n-2},0,z_{n})\in\partial D,$$ $$d_D(z)\leq\min_{j=n-1,n}|z_{j}|.$$ 

Case 2.2. All the points $$(z',0,\ldots,0,\underline{z_p},0,\ldots,0,\underline{z_q},0,\ldots,0),\,k+1\leq p<q\leq n$$ belong to $\partial D$. We see, by induction, that in the $s$-th step ($s=3,\ldots,n-k-1$) all points 
$$(z',0,\ldots,0,\underline{z_{p_1}},0,\ldots,0,\underline{z_{p_s}},0,\ldots,0),\,k+1\leq p_1<\ldots p_s\leq n,$$ except possibly one, belong to $\partial D$. 

If one of these points, say $(z',0,\ldots,0,z_{n-s+1},\ldots,z_n)$, does not belong to $\partial D$ then it belongs to $D$ and $$P(\zeta_0,r)\cap(\mathbb{C}^{n-s}\times\mathbb{C}^s_*)\subset D,$$ $$k_D(z_0,z)\leq\frac{1}{2}\log(-\log\min_{j=n-s+1,\ldots,n}|z_j|)+C,$$ 
$$(z',z_{k+1},\ldots,z_{n-s},z_{n-s+1},\ldots,z_{j-1},0,z_{j+1},\ldots,z_n)\in\partial D,\,j=n-s+1,\ldots,n,$$
$$d_D(z)\leq\min_{j=n-s+1,\ldots,n}|z_{j}|,$$ which finishes the proof in the case $s.1$. 

If all the points $$(z',0,\ldots,0,\underline{z_{p_1}},0,\ldots,0,\underline{z_{p_s}},0,\ldots,0),\,k+1\leq p_1<\ldots p_s\leq n$$ belong to $\partial D$ then we ``jump'' from the case $s$.2 to the case $(s+1).1$, getting finally in the case $(n-k-1).1$ $$(z',0,z_{k+2},\ldots,z_n)\in D,$$ $$P(\zeta_0,r)\cap(\mathbb{C}^{k+1}\times\mathbb{C}^{n-k-1}_*)\subset D,$$ $$k_D(z_0,z)\leq\frac{1}{2}\log(-\log\min_{j=k+2,\ldots,n}|z_j|)+C,$$ 
$$(z',z_{k+1},z_{k+2},\ldots,z_{j-1},0,z_{j+1},\ldots,z_n)\in\partial D,\,j=k+2,\ldots,n,$$
$$d_D(z)\leq\min_{j=k+2,\ldots,n}|z_{j}|$$ or in the case $(n-k-1).2$ $$(z',z_{k+1},\ldots,z_{j-1},0,z_{j+1},\ldots,z_{n})\in\partial D,\,j=k+1,\ldots,n.$$ This property let us estimate $d_D(z)$ from above by $\min_{j=k+1,\ldots,n}|z_j|$ and use \eqref{11} to finish the proof.
\end{proof}
\begin{proof}[\sc Proof of Theorem \ref{3}] 
The proof has two main parts; in the first the claim is proved for $\zeta_0\in\partial D\cap\mathbb{C}^{n}_{*}$ thanks to the effective formulas for the Kobayashi distance in special domains and in the second part the remaining case is amounted to the lower-dimensional situation with a boundary point having all non-zero coordinates.

Let $\zeta_0\in\partial D\cap\mathbb{C}^{n}_{*}$ and consider $z\in D\cap\mathbb{C}^{n}_{*}$ close to $\zeta_0$. From the convexity of the set $\log{D}$ there exist $\alpha\in\mathbb{R}^{n}$ and $c>0$ such that the hyperplane
$$\{x\in\mathbb{R}^{n}:\langle\alpha,x\rangle_{\mathbb{R}^{n}}=\log c\}$$ contains point $\log|\zeta_0|$ and $\log{D}$ lies on the one side of this hyperplane. Assume, without loss of generality, that this side is $\{x\in\mathbb{R}^{n}:\langle\alpha,x\rangle_{\mathbb{R}^{n}}<\log c\}$ since in the case of $\log{D}\subset\{x\in\mathbb{R}^{n}:\langle\alpha ',x\rangle_{\mathbb{R}^{n}}>\log c'\}$ it suffices to define $$\alpha:=-\alpha '\text{ and }c:=1/c'.$$
Therefore $$\{( e^{x_{1}},\ldots,e^{x_{n}}):x\in \log D\}\subset\{w\in\mathbb{C}^{n}:|w|^{\alpha}<c\}=:D_{\alpha,c}\footnote{These sets are called \textit{elementary Reinhardt domains}.},$$ 
where by a point satisfying the condition $|w|^{\alpha}<c$ we mean such point $w$ whose coordinate $w_{j}$ is non-zero when $\alpha_j<0$ (and satisfies $|w|^{\alpha}<c$ in the usual sense). To affirm that $D\subset D_{\alpha,c}$, we have to check that the above restriction for points $w$ does not remove from $D$ points with some zero coordinates. Indeed, if there is no such inclusion, we can assume that the order of zero coordinates of point $w\in{D}$ and negative terms of the sequence $\alpha$ is as follows: $$w_{1},\ldots,w_{k}\neq 0,\, w_{k+1},\ldots,w_{n}=0$$ $$\alpha_{k+1},\ldots,\alpha_{l}\geq 0,\,\alpha_{l+1},\ldots,\alpha_{n}<0,$$ 
where $1\leq k\leq l<n$. In some neighbourhood of the point $w$ contained in ${D}$ there exist points $v\in\mathbb{C}_{*}^{n}$ with coordinates $v_{j}$ such that $$|v_{1}|,\ldots,|v_{l}|>\varepsilon>0$$ and $|v_{l+1}|,\ldots,|v_{n}|$ arbitrarily close to zero (i.e. moved from $w$ in a direction of subspace $\{0\}^{l}\times\mathbb{C}^{n-l}$ and next moved from it by a constant vector in the direction $\mathbb{C}^{l}\times\{0\}^{n-l}$). Then there exist points $u\in\log{D}$ whose coordinates $u_{j}$ satisfy $$u_{1},\ldots,u_{l}>\log\varepsilon>-\infty$$
however $u_{l+1},\ldots,u_{n}$ are arbitrarily close to $-\infty$. But it contradicts a fact that values of the expression
$$\sum_{j=l+1}^{n}\alpha_{j}u_{j}$$ are for these points $u$ bounded from above by a constant $\log c-\sum_{j=1}^{l}\alpha_{j}\log\varepsilon$. 

We will use effective formulas for the Kobayashi distance in domains $D_{\alpha,c}$ \cite{zw1}. Define $$l:=\#\{j=1,\ldots,n:\alpha_{j}<0\}$$ and $$\widetilde{\alpha}:=\min\{\alpha_{j}:\alpha_{j}>0\}\,\textnormal{ if }\,l<n.$$ We first consider a situation $l<n$. The formula in this case gives $$k_{D}(z_{0},z)\geq k_{D_{\alpha,c}}(z_{0},z)\geq p\left(0,\frac{|z|^{\alpha/\widetilde{\alpha}}}{c^{1/\widetilde{\alpha}}}\right)+C.$$ But $$z=\zeta_{D_{\alpha,c}}(z)-d_{D_{\alpha,c}}(z)\nu_{D_{\alpha,c}}(z)$$ and hence $$|z|^{\alpha/\widetilde{\alpha}}=\prod\limits_{j=1}^{n}|\zeta_{D_{\alpha,c}}(z)_{j}-d_{D_{\alpha,c}}{(z)}
\nu_{D_{\alpha,c}}(z)_{j}|^{\alpha_{j}/\widetilde{\alpha}}=c^{1/\widetilde{\alpha}}-\rho(z)d_{D_{\alpha,c}}(z)$$ for some bounded positive function $\rho$. Thus $$p\left(0,\frac{|z|^{\alpha/\widetilde{\alpha}}}{c^{1/\widetilde{\alpha}}}\right)=
p\left(0,1-\frac{\rho(z)}{c^{1/\widetilde{\alpha}}}d_{D_{\alpha,c}}(z)\right)\geq
-\frac{1}{2}\log\left(\frac{\rho(z)}{c^{1/\widetilde{\alpha}}}d_{D_{\alpha,c}}(z)\right)\geq$$$$\geq
-\frac{1}{2}\log d_{D_{\alpha,c}}(z)+C.$$
We will show that $$d_{D_{\alpha,c}}(z)\approx d_{{D}}(z)\textnormal{ as }z\to\zeta_{0}
\textnormal{ non-tangentially}.$$
By the definition there exists a cone $\mathcal{A}$ with a vertex $\zeta_{0}$ and a semi-axis $-\nu_{D_{\alpha,c}}(\zeta_{0})$ which contains considered points $z$. Thanks to the $\mathcal{C}^{1}$-smoothness of ${D}$ we have a cone $\mathcal{B}$ with the vertex $\zeta_{0}$ and the semi-axis $-\nu_{D_{\alpha,c}}(\zeta_{0})$, whose intersection with some neighbourhood of the point $\zeta_{0}$ is contained in ${D}$ and contains in its interior the cone $\mathcal{A}$. Therefore
$$1\geq\frac{{d_{{D}}(z)}}{d_{D_{\alpha,c}}(z)}=\frac{\|z-\zeta_{{D}}(z)\|}{\|z-\zeta_{D_{\alpha,c}}(z)\|}\geq
\frac{\|z- \zeta_{{D}}(z)\|}{\|z-\zeta_{0}\|}\geq\frac{\|z-\zeta_{\mathcal{B}}(z)\|}{\|z-\zeta_{0}\|}=$$ $$=\sin\angle(z,\zeta_{0},\zeta_{\mathcal{B}}(z))\geq\sin\theta,$$ where $\angle(X,Y,Z)$ is an angle with vertex $Y$, whose arms contain points $X$, $Z$ and $\theta$ is an angle between these generatrices of cones $\mathcal{A}$, $\mathcal{B}$ which lie in one plane with the axis of both cones.\footnote{In other words, $\theta$ is a difference of angles appearing in the definitions of the cones $\mathcal{B}$, $\mathcal{A}$.}

The second case $l=n$ gives $$k_{D}(z_{0},z)\geq k_{D_{\alpha,c}}(z_{0},z)\geq p\left(0,\frac{|z|^{\alpha}}{c}\right)+C.$$ Similarly as before $$|z|^{\alpha}=\prod\limits_{j=1}^{n}|\zeta_{D_{\alpha,c}}(z)_{j}-d_{D_{\alpha,c}}{(z)}
\nu_{D_{\alpha,c}}(z)_{j}|^{\alpha_{j}}=c-\sigma(z)d_{D_{\alpha,c}}(z)$$ with a bounded positive function $\sigma$. Hence $$p\left(0,\frac{|z|^{\alpha}}{c}\right)\geq-\frac{1}{2}\log d_{D_{\alpha,c}}(z)+C\geq-\frac{1}{2}\log d_D(z)+C.$$ 

Now, take $\zeta_{0}\in\partial D\setminus\mathbb{C}_{*}^{n}$. We may assume that the first $k$ coordinates of $\zeta_0$ are non-zero and the last $n-k$ are zero, where $0\leq k\leq n-1$. Notice that $k\neq 0$. Indeed, the assumption $k=0$ is equivalent to $0\in\partial{D}$. 
Using Facts \ref{4} and \ref{5} we see that the $\mathcal{C}^{1}$-smoothness of $D$ implies, thanks to the Fu condition, ${D}\cap V_{j}^{n}\neq\emptyset$ for $j=1,\ldots,n$. Hence $D$ is complete i.e. $0\in D$ --- a contradiction. Finally, point $\zeta_{0}$ has a form $$\zeta_{0}=(\zeta_{01},\ldots,\zeta_{0k},0,\ldots,0),\,\zeta_{0j}\neq 0,\,1\leq j\leq k\leq n-1.$$ Consider the projection $\pi_{k}:\mathbb{C}^{n}\longrightarrow\mathbb{C}^{k}$ i.e. $$\pi_{k}(z)=(z_{1},\ldots,z_{k}).$$
We will show that ${D}_{k}:=\pi_{k}({D})$ is a $\mathcal{C}^{1}$-smooth pseudoconvex Reinhardt domain. A Reinhardt property is clear for ${D}_{k}$. To affirm the pseudoconvexity of $D_k$ it suffices to show that $${D}_{k}\times\{0\}^{n-k}={D}\cap(\mathbb{C}^{k}\times\{0\}^{n-k}).$$ Inclusion $${D}_{k}\times\{0\}^{n-k}\supset{D}\cap(\mathbb{C}^{k}\times\{0\}^{n-k})$$
is obvious. To prove the opposite inclusion we will use Facts \ref{4} and \ref{5} again. We have ${D}\cap V_{j}^{n}\neq\emptyset$ for $j=k+1,\ldots,n$, so $D$ is complete in $j$-th direction for $j=k+1,\ldots,n$. Take some $z\in D_{k}\times\{0\}^{n-k}$. Then $z=(z_{1},\ldots,z_{k},0\ldots,0)$ and $(z_{1},\ldots,z_{k},\widetilde{z}_{k+1},\ldots,\widetilde{z}_{n})\in D$ for some $\widetilde{z}_{k+1},\ldots,\widetilde{z}_{n}\in\mathbb{C}$. Thus $(z_{1},\ldots,z_{k},0\ldots,0)\in D$ i.e. $z\in D\cap(\mathbb{C}^{k}\times\{0\}^{n-k})$. 

The local defining function for ${D}_{k}$ at point $\zeta\in\partial{D}_{k}$ is $$\widetilde{\rho}(z_{1},\ldots,z_{k}):=\rho(z_{1},\ldots,z_{k},0,\ldots,0),\,(z_{1},\ldots,z_{k})\in\pi_{k}(U)\cap{D}_{k},$$
where $\rho:U\longrightarrow\mathbb{R}$ is the local defining function for ${D}$ at point $(\zeta,0,\ldots,0)$.
Indeed, $\nabla\widetilde{\rho}\neq 0$ since
\begin{itemize}
\item $\nabla\rho\neq 0$; 
\item $\frac{\partial\widetilde{\rho}}{\partial\overline{z}_{j}}=
\frac{\partial\rho}{\partial\overline{z}_{j}}$ for $j=1,\ldots,k$;
\item $\frac{\partial\rho}{\partial\overline{z}_{j}}=0$ for $j=k+1,\ldots,n$,
\end{itemize}
however the two remaining conditions for a defining function follow easy from the definition of $\widetilde{\rho}$. 

If $z$ tends to $\zeta_{0}$ non-tangentially in a cone $\mathcal{A}\subset\mathbb{C}^{n}$ then $\pi_{k}(z)$ tends to $\pi_{k}(\zeta_{0})\in\mathbb{C}_{*}^{k}$ non-tangentially in a cone $\pi_{k}(\mathcal{A})\subset\mathbb{C}^{k}$. From the previous case $$k_{{D}}(z_{0},z)\geq k_{{D}_{k}}(\pi_{k}(z_{0}),\pi_{k}(z))\geq-\frac{1}{2}\log d_{{D}_{k}}(\pi_{k}(z))+C.$$
Hence to finish the proof it suffices to show that $$d_{{D}_{k}}(\pi_{k}(z))\lesssim d_{{D}}(z).$$
Consider a cone $\mathcal{B}$ with vertex $\zeta_{0}$ and semi-axis $-\nu_{D_{\alpha}}(\zeta_{0})$ whose intersection with some neighbourhood of the point $\zeta_{0}$ is contained in ${D}$ and contains in its interior the cone $\mathcal{A}$. Then $$1\geq\frac{d_{\mathcal{B}}(z)}{d_{{D}}(z)}=\frac{\|z-\zeta_{\mathcal{B}}(z)\|}{\|z-\zeta_{{D}}(z)\|}\geq
\frac{\|z-\zeta_{\mathcal{B}}(z)\|}{\|z-\zeta_{0}\|}=\sin\angle(z,\zeta_{0},\zeta_{\mathcal{B}}(z))\geq\sin\theta,$$ where $\theta$ is an angle between these generatrices of the cones $\mathcal{A}$, $\mathcal{B}$ which lie in one plane with the axis of both cones. Analogously
$$1\geq\frac{d_{\pi_{k}(\mathcal{B})}(\pi_{k}(z))}{d_{{D}_{k}}(\pi_{k}(z))}\geq\sin\theta',$$
where $\theta'$ depends only on $\mathcal{B}$. 
Therefore $$\frac{d_{{D}_{k}}(\pi_{k}(z))}{d_{{D}}(z)}\approx\frac{d_{\pi_{k}(\mathcal{B})}(\pi_{k}(z))}{d_{\mathcal{B}}(z)},$$
however $$\frac{d_{\pi_{k}(\mathcal{B})}(\pi_{k}(z))}{d_{\mathcal{B}}(z)}=
\frac{\|\pi_{k}(z)-\zeta_{\pi_{k}(\mathcal{B})}(\pi_{k}(z))\|}{\|z-\zeta_{\mathcal{B}}(z)\|}=$$\medskip
$$=\frac{\|\pi_{k}(z)-\pi_{k}(\zeta_{0})\|\sin\angle(\pi_{k}(z),\pi_{k}(\zeta_{0}),\zeta_{\pi_{k}(\mathcal{B})}(\pi_{k}(z)))}
{\|z-\zeta_{0}\|\sin\angle(z,\zeta_{0},\zeta_{\mathcal{B}}(z))}\leq$$\medskip
$$\leq\frac{\|\pi_{k}(z)-\pi_{k}(\zeta_{0})\|}{\|z-\zeta_{0}\|\sin\theta}\leq
\frac{1}{\sin\theta}.\eqno\qedhere$$
\end{proof}
\begin{rem}
The estimate from below by $-\frac{1}{2}\log d_D+C$ for the Carath\'eodory (pseudo) distance $c_D$ is not true even for smooth bounded complete pseudoconvex Reinhardt domain $D$ and its boundary point $\zeta_0\in\mathbb{C}^{n}_{*}$.
\end{rem}
\begin{proof}[\sc Proof.]
Consider a domain
$$D:=\left\{(z_1,z_2)\in\mathbb{C}^2:|z_1|<R_1,\,|z_2|<R_2,\,|z_1||z_2|^\alpha<R_3\right\},$$ where $R_1,R_2,R_3>0$, $\alpha\in(\mathbb{R}\setminus\mathbb{Q})_+$ and $R_{1}R_{2}^{\alpha}>R_3$. Fix $\zeta_0\in\partial D$ such that $|\zeta_{01}|<R_1,|\zeta_{02}|<R_2$. This domain is not smooth. Since $$\log D=\left\{(x_1,x_2)\in\mathbb{R}^2:x_1<\log R_1,\,x_2<\log R_2,\,x_1+\alpha x_2<\log R_3\right\},$$ it is easy to construct smooth bounded convex domain $E\subset\mathbb{R}^2$ such that $\log D\subset E$ and $\partial E$ contains the skew segment $$(\partial\log D)\cap\{(x_1,x_2)\in\mathbb{R}^2:x_1+\alpha x_2=\log R_3\}.$$ Let $\widetilde{D}\subset\mathbb{C}^2$ be a complete Reinhardt domain such that $\log\widetilde{D}=E$. Then $\widetilde{D}$ is bounded, smooth and, thanks to Fact \ref{4}, pseudoconvex. Moreover, $D\subset\widetilde{D}$ and $D$, $\widetilde{D}$ are identic in the neighbourhood of their common boundary point $\zeta_0$.

We have from the Proposition 4.3.2 in \cite{zw1} $$a_\lambda:=\frac{g_D(\lambda\zeta_0,0)}{\log|\lambda|}\to\infty\text{ as }\lambda\to\partial\mathbb{D},$$ where $g_D$ is the pluricomplex Green function (general properties of the Carath\'eodory (pseudo)distance and the pluricomplex Green function one can find e.g. in \cite{jp} and \cite{zw1}). Certainly $$d_{\widetilde{D}}(\lambda\zeta_0)=d_D(\lambda\zeta_0)\approx 1-|\lambda| \text{ as } |\lambda|\to 1$$ and $$c_D(\lambda\zeta_0,0)\leq\tanh^{-1}\exp g_D(\lambda\zeta_0,0).$$ Therefore, if there exists a constant $C>0$ such that $$c_{\widetilde{D}}(\lambda\zeta_0,0)\geq-\frac{1}{2}\log d_{\widetilde{D}}(\lambda\zeta_0)+C,\,|\lambda|\to 1$$ then for $|\lambda|\to 1$ $$c_{D}(\lambda\zeta_0,0)\geq-\frac{1}{2}\log d_{D}(\lambda\zeta_0)+C,$$ $$-\frac{1}{2}\log(1-|\lambda|)+C\leq\tanh^{-1}|\lambda|^{a_\lambda},$$ $$\frac{1}{1-|\lambda|}\leq\frac{C'}{1-|\lambda|^{a_\lambda}}$$ with a constant $C'>0$. For $|\lambda|$ sufficiently close to 1 we have $a_\lambda\geq C'+1$. Therefore $$\frac{1}{1-|\lambda|}\leq\frac{C'}{1-|\lambda|^{C'+1}}$$ or equivalently $$\frac{1-|\lambda|^{C'+1}}{1-|\lambda|}\leq C'.$$ The left-hand side tends to $C'+1$ as $|\lambda|\to 1$.
\end{proof}
\section{Open problems}
\begin{enumerate}[(1)]
\item Can we improve the estimate from Theorem \ref{1} to $-\frac{1}{2}\log d_{{D}}(z)+C$?
\item Let ${D}\subset\mathbb{C}^{n}$ be a pseudoconvex Reinhardt domain and $\zeta_{0}\in\partial D\cap\mathbb{C}^{n}_{*}$. Does it implies that for some constant $C$ the inequality $$k_{{D}}(z_{0},z)\geq-\frac{1}{2}\log d_{{D}}(z)+C$$ holds if $z\in D$ tends to $\zeta_{0}$?
\item Is it true for pseudoconvex Reinhardt domains $D\subset\mathbb{C}^{n}$ that if $$\#\{j:\zeta_{0j}=0\text{ and } D\cap V_{j}^{n}=\emptyset\}=0$$
then $$k_{D}(z_{0},z)\geq-\frac{1}{2}\log d_{D}(z)+C$$ and in the opposite case
$$k_{D}(z_{0},z)\geq\frac{1}{2}\log(-\log d_{D}(z))+C$$ for $z\in D$ near $\zeta_0\in\partial D$?\\
\end{enumerate}

\smallskip
\textsc{Acknowledgements.} I would like to thank Professor Włodzimierz Zwonek for introduction to the problem and helpful suggestions.

\end{document}